\documentclass[11pt]{article}

\usepackage{latexsym}
\usepackage{amssymb}

\newtheorem{thm}{Theorem}

\newtheorem{dfn}{Definition}

\begin{document}
{
\begin{center}
{\Large\bf
Orthogonal polynomials related to some Jacobi-type pencils.}
\end{center}
\begin{center}
{\bf S.M. Zagorodnyuk}
\end{center}

\section{Introduction.}

The theory of orthogonal polynomials on the real line is a classical subject which has a huge ammount of
contributions and applications in various disciplines~\cite{cit_50000_Gabor_Szego},\cite{cit_20000_Suetin},\cite{cit_3000_Chihara}. 
Nowadays it attracts a lot of researches who use
various new techniques~\cite{cit_5000_Ismail},\cite{cit_10000_Simon}. In this paper we shall study a generalization of the class of 
orthogonal polynomials on the real line.
For this aim we shall need the following basic definition.

\begin{dfn}
\label{d1_1}
A set $\Theta = \left(
J_3, J_5, \alpha, \beta
\right)$,
where $\alpha>0$, $\beta\in\mathbb{R}$, $J_3$ is a Jacobi matrix and
$J_5$ is a semi-infinite real symmetric five-diagonal matrix with positive numbers on the second subdiagonal,
is said to be
\textbf{a Jacobi-type pencil (of matrices)}.
\end{dfn}
As it follows from this definition the matrices $J_3$ and $J_5$ have the following form:
\begin{equation}
\label{f1_5}
J_3 =
\left(
\begin{array}{cccccc}
b_0 & a_0 & 0 & 0 & 0 & \cdots\\
a_0 & b_1 & a_1 & 0 & 0 & \cdots\\
0 & a_1 & b_2 & a_2 & 0 &\cdots\\ 
\vdots & \vdots & \vdots & \ddots \end{array}
\right),\qquad a_k>0,\ b_k\in\mathbb{R},\ k\in\mathbb{Z}_+;
\end{equation}

\begin{equation}
\label{f1_10}
J_5 =
\left(
\begin{array}{ccccccc}
\alpha_0 & \beta_0 & \gamma_0 & 0 & 0 & 0 & \cdots\\
\beta_0 & \alpha_1 & \beta_1 & \gamma_1 & 0 & 0 & \cdots\\
\gamma_0 & \beta_1 & \alpha_2 & \beta_2 & \gamma_2 & 0 & \cdots\\ 
0 & \gamma_1 & \beta_2 & \alpha_3 & \beta_3 & \gamma_3 & \cdots \\
\vdots & \vdots & \vdots &\vdots & \ddots \end{array}
\right),\ \alpha_n,\beta_n\in\mathbb{R},\ \gamma_n>0,\ n\in\mathbb{Z}_+.
\end{equation}
Jacobi matrices and the corresponding operators are closely related to orthogonal polynomials on the real line
and are intensively studied~\cite{cit_70000_T}.
Five-diagonal matrices may be viewed as $(2\times 2)$ block Jacobi matrices. These matrices are related to
$(2\times 2)$ matrix orthogonal polynomials on the real line~\cite{cit_4000_D_P_S}, as well as to orthogonal polynomials on
radial rays in the complex plane (see, e.g.~\cite{cit_90000_Z},\cite{cit_3500_Ch_Z} and references therein).

With a Jacobi-type pencil of matrices $\Theta$ we shall associate a system of polynomials
$\{ p_n(\lambda) \}_{n=0}^\infty$, such that
\begin{equation}
\label{f1_15}
p_0(\lambda) = 1,\quad p_1(\lambda) = \alpha\lambda + \beta,
\end{equation}
and
\begin{equation}
\label{f1_20}
(J_5 - \lambda J_3) \vec p(\lambda) = 0,
\end{equation}
where $\vec p(\lambda) = (p_0(\lambda), p_1(\lambda), p_2(\lambda),\cdots)^T$. Here the superscript $T$ means the transposition.
Polynomials $\{ p_n(\lambda) \}_{n=0}^\infty$ are said to be \textbf{associated to the Jacobi-type pencil of matrices $\Theta$}.

By choosing all possible Jacobi-type pencils of matrices we obtain a class $\mathfrak{K}$ which consists of the associated
systems of polynomials. The class $\mathfrak{K}$ contains the class $\mathfrak{R}$ of all systems of orthonormal polynomials on the real line 
with $p_0=1$ (and positive leading coefficients).
In fact, for each system of orthonormal polynomials on the real line with $p_0=1$ one may choose $J_3$ to be the corresponding Jacobi matrix
(which elements are the recurrence coefficients),  $J_5 = J_3^2$, and $\alpha,\beta$  being the coefficients of $p_1$
($p_1(\lambda) = \alpha\lambda + \beta$).

In the case of all bounded coefficients of the matrices $J_3$, $J_5$, these matrices define, in the usual manner, bounded operators
on the space $l_2$. These operators will be denoted by the same letters as matrices.
In this case, the expression $J_5 - \lambda J_3$ form the linear operator pencil. For the theory of polynomial operator pencils
we refer to the book~\cite{cit_7000_Markus}.

In the general case, the matrices $J_3$, $J_5$ allow us to define operators $J_{3,0}$, $J_{5,0}$ on the set of all
finite vectors from $l_2$ (i.e. complex vectors with all but finite number coefficients zeros).
These operators will be used later to define an operator of the pencil and to study its properties.

Relation~(\ref{f1_20}) may be written in the following scalar form:
$$ \gamma_{n-2} p_{n-2}(\lambda) + (\beta_{n-1}-\lambda a_{n-1}) p_{n-1}(\lambda) + (\alpha_n-\lambda b_n) p_n(\lambda) +
$$
\begin{equation}
\label{f1_30}
+ (\beta_n-\lambda a_n) p_{n+1}(\lambda) + \gamma_n p_{n+2}(\lambda) = 0,\qquad n\in\mathbb{Z}_+,
\end{equation}
where $p_{-2}(\lambda) = p_{-1}(\lambda) = 0$, $\gamma_{-2} = \gamma_{-1} = \alpha_{-1} = \beta_{-1} = 0$.
Recurrent relation~(\ref{f1_30}) with the initial conditions~(\ref{f1_15}) uniquely determine the associated polynomials
of any Jacobi-type pencil of matrices. Moreover, the polynomial $p_n$ has degree $n$ and a positive leading coefficient ($n\in\mathbb{Z}_+$).
On the other hand, it is clear that the associated polynomials do not determine the pencil. For example, multiplying
$J_3$ and $J_5$ by a positive constant does not change the associated polynomials.

Our first aim is to show that the inclusion of $\mathfrak{R}$ into $\mathfrak{K}$ is strict, i.e. $\mathfrak{R}\not=\mathfrak{K}$.
For this purpose we construct our basic example of polynomials associated to a pencil. Moreover, the associated polynomials
admit an explicit representation.
Our second aim is to obtain some  orthonormality relations for the associated polynomials of an arbitrary Jacobi-type pencil.
For that purpose we shall introduce an operator of the pencil.

{\bf Notations. }
As usual, we denote by $\mathbb{R}, \mathbb{C}, \mathbb{N}, \mathbb{Z}, \mathbb{Z}_+$,
the sets of real numbers, complex numbers, positive integers, integers and non-negative integers,
respectively. By $\mathbb{P}$ we denote the set of all polynomials with complex coefficients.

\noindent
By $l_2$ we denote the usual Hilbert space of all complex sequences $c = (c_n)_{n=0}^\infty = (c_0,c_1,c_2,...)^T$ with the finite norm
$\| c \|_{l_2} = \sqrt{\sum_{n=0}^\infty |c_n|^2}$. Here $T$ means the transposition.
The scalar product of two sequences $c = (c_n)_{n=0}^\infty, d = (d_n)_{n=0}^\infty\in l_2$ is given by
$(c,d)_{l_2} = \sum_{n=0}^\infty c_n \overline{d_n}$. 
We denote $\vec e_m = (\delta_{n,m})_{n=0}^\infty\in l_2$, $m\in\mathbb{Z}_+$.
By $l_{2,fin}$ we denote the set of all finite vectors from $l_2$, i.e. vectors with all but finite number components zeros.

\noindent
By $\mathfrak{B}(\mathbb{R})$ we denote the set of all Borel subsets of $\mathbb{R}$.
If $\sigma$ is a (non-negative) bounded measure on $\mathfrak{B}(\mathbb{R})$ then by $L^2_\sigma$ we denote a Hilbert space
of all (classes of equivalences) of complex-valued functions $f$ on $\mathbb{R}$ with a finite norm
$\| f \|_{L^2_\sigma} = \sqrt{ \int_\mathbb{R} |f(x)|^2 d\sigma }$.
The scalar product of two functions $f,g\in L^2_\sigma$ is given by
$(f,g)_{L^2_\sigma} = \int_{\mathbb{R}} f(x) \overline{g(x)} d\sigma$. 
By $[ f ]$  we denote the class of equivalence in $L^2_\sigma$ which contains the representative $f$.

\noindent
If H is a Hilbert space then $(\cdot,\cdot)_H$ and $\| \cdot \|_H$ mean
the scalar product and the norm in $H$, respectively.
Indices may be omitted in obvious cases.
For a linear operator $A$ in $H$, we denote by $D(A)$
its  domain, by $R(A)$ its range, by $\mathop{\rm Ker}\nolimits A$
its null subspace (kernel), and $A^*$ means the adjoint operator
if it exists. If $A$ is invertible then $A^{-1}$ means its
inverse. $\overline{A}$ means the closure of the operator, if the
operator is closable. If $A$ is bounded then $\| A \|$ denotes its
norm.
For a set $M\subseteq H$
we denote by $\overline{M}$ the closure of $M$ in the norm of $H$.
By $\mathop{\rm Lin}\nolimits M$ we mean
the set of all linear combinations of elements of $M$,
and $\mathop{\rm \overline{span}}\nolimits M := \overline{ \mathop{\rm Lin}\nolimits M }$.
By $E_H$ we denote the identity operator in $H$, i.e. $E_H x = x$,
$x\in H$. If $H_1$ is a subspace of $H$, then $P_{H_1} =
P_{H_1}^{H}$ is an operator of the orthogonal projection on $H_1$
in $H$.

\section{The basic example.}

Denote by $\Theta_1$ the Jacobi-type pencil with $\alpha = \beta = \sqrt{2}$;
$a_k = \sqrt{2}$, $b_k = 2$, $k\in\mathbb{Z}_+$;
$\alpha_n = \beta_n = 0$, $\gamma_n = 1$, $n\in\mathbb{Z}_+$.
Recurrent relation~(\ref{f1_30}) for the associated polynomials $\{ p_n(\lambda) \}_{n=0}^\infty$ takes the following
form:
\begin{equation}
\label{f2_10}
p_{n-2}(\lambda) - \sqrt{2} \lambda p_{n-1}(\lambda) - 2 \lambda p_n(\lambda) 
- \sqrt{2} \lambda  p_{n+1}(\lambda) + p_{n+2}(\lambda) = 0,\qquad n\in\mathbb{Z}_+,
\end{equation}
with the initial conditions:
\begin{equation}
\label{f2_15}
p_0(\lambda) = 1,\quad p_1(\lambda) = \sqrt{2}\lambda + \sqrt{2}.
\end{equation}

\begin{thm}
\label{t2_1}
The associated polynomials $\{ p_n(\lambda) \}_{n=0}^\infty$ of the pencil $\Theta_1$ have the following
representation:
$$ p_n(\sqrt{2}t-1) = T_n(t) + t U_{n-1}(t) - \frac{1}{2} 
\frac{ U_{n-1}(t) - U_{n-1}\left( -\frac{1}{ \sqrt{2} } \right) }
{ t + \frac{1}{\sqrt{2}} }, $$
\begin{equation}
\label{f2_20}
n\in\mathbb{Z}_+,\ t\in \left(-1, -\frac{1}{\sqrt{2}} \right) \cup \left( -\frac{1}{\sqrt{2}}, 1 \right).
\end{equation}
Here $T_n(t) = \cos(n\arccos t)$, $U_n(t) = \frac{ \sin((n+1)\arccos t) }{ \sqrt{1-t^2} }$
are Chebyshev's polynomials of the first and the second kind, respectively ($U_{-1}=0$).
\end{thm}
\textbf{Proof.}
As usual, we shall seek for the solution $\{ y_n \}_{n=0}^\infty$ of the difference equation
(related to~(\ref{f2_10}))
\begin{equation}
\label{f2_23}
y_{n-2} - \sqrt{2} \lambda y_{n-1} - 2 \lambda y_n 
- \sqrt{2} \lambda  y_{n+1} + y_{n+2} = 0,\qquad n = 2,3,...,
\end{equation}
in the following form:
$y_n = w^n$, $n\in\mathbb{Z}_+$, $w=w(\lambda)\in\mathbb{C}$. 
Here $\lambda\in\mathbb{C}$ is a fixed parameter.
We obtain the following characteristic equation:
\begin{equation}
\label{f2_25}
w^4 - \sqrt{2} \lambda w^3 - 2\lambda w^2 - \sqrt{2} \lambda w + 1 = 0.
\end{equation}
We may factorize the expression on the left to obtain that
\begin{equation}
\label{f2_30}
(w^2 + \sqrt{2} w + 1) (w^2 - \sqrt{2} (\lambda + 1) w + 1) = 0.
\end{equation}
Thus, we have the following roots:
$$ w_{1,2} = \frac{ \sqrt{2} }{2} (-1\mp i),\quad
w_{3,4} = \frac{ \sqrt{2} }{2} \left(
\lambda+1 \pm \sqrt{\lambda^2 + 2\lambda -1}
\right). $$
Here and in what follows, for each complex $\lambda$ we fix an \textit{arbitrary} value of the square root $\sqrt{\lambda^2 + 2\lambda -1}$.
We do not assume that these values form some branch or require other conditions. 
Set
$$ r_n(\lambda) = C_1(\lambda) w_1^n + C_2(\lambda) w_2^n + C_3(\lambda) w_3^n + C_4(\lambda) w_4^n,\qquad n\in \mathbb{Z}_+, $$  
where $C_j(\lambda)$ are arbitrary complex-valued functions of $\lambda$.
The functions $r_n(\lambda)$ satisfy relation~(\ref{f2_10}) for $n=2,3,...$.
Moreover,  $r_n(\lambda)$ satisfy relation~(\ref{f2_10}) with $n=0,1$ and the initial conditions~(\ref{f2_15}) if and only if 
the coefficients $C_n(\lambda)$
satisfy the following linear system of equations: 
\begin{equation}
\label{f2_35}
\left\{
\begin{array}{ccccccc}
(-\lambda + (\lambda+1)i) C_1 + (-\lambda - (\lambda+1)i) C_2 +
(-\lambda + \sqrt{\lambda^2 + 2\lambda -1} ) C_3 +  \\
+ (-\lambda - \sqrt{\lambda^2 + 2\lambda -1} ) C_4 = 0, \\
(1-i) C_1 + (1+i) C_2 +
(-\lambda -1 + \sqrt{\lambda^2 + 2\lambda -1} ) C_3 + \\
+ (-\lambda - 1 - \sqrt{\lambda^2 + 2\lambda -1} ) C_4 = 0, \\
C_1 + C_2 + C_3 + C_4 = 1, \\
(-1-i) C_1 + (-1+i) C_2 +
(\lambda + 1 + \sqrt{\lambda^2 + 2\lambda -1} ) C_3 + \\
+ (\lambda + 1 - \sqrt{\lambda^2 + 2\lambda -1} ) C_4 = 2\lambda + 2.\end{array}
\right.
\end{equation}
The determinant $\Delta$ of this system is equal to
$-8 (\lambda + 2)^2 \sqrt{\lambda^2 + 2\lambda -1} i$.
Thus, this linear system has a solution if $\lambda\not= -2, -1\pm\sqrt{2}$.
Then
$$ C_{1,2} = \mp\frac{ 1 }{2(\lambda + 2)i},\
C_{3,4} = \frac{1}{2} \pm \frac{ \lambda^2 + 3\lambda + 1 }{2(\lambda + 2)\sqrt{\lambda^2 + 2\lambda -1}}, $$
and we come to the following representation of the associated polynomials:
$$ p_n(\lambda) = 
\frac{1}{ \lambda + 2 } \sin\left( \frac{3\pi}{4} n \right) +
2^{-\frac{n}{2} - 1}
\left(
(\lambda + 1 + \sqrt{\lambda^2 + 2\lambda -1})^n
+
\right.
$$
$$
+ (\lambda + 1 - \sqrt{\lambda^2 + 2\lambda -1})^n + \frac{\lambda^2 + 3\lambda + 1}{ (\lambda + 2) \sqrt{\lambda^2 + 2\lambda -1} } * $$
$$ \left. * \left(
(\lambda + 1 + \sqrt{\lambda^2 + 2\lambda -1})^n
-
(\lambda + 1 - \sqrt{\lambda^2 + 2\lambda -1})^n
\right)
\right), $$
\begin{equation}
\label{f2_40}
n\in\mathbb{Z}_+,\ \lambda\in\mathbb{C}\backslash\{ -2, -1\pm\sqrt{2} \}.
\end{equation}
In what follows we suppose that $\lambda\in(-1-\sqrt{2}, -2)\cup (-2, -1+\sqrt{2})$. Then 
$t := \frac{\lambda + 1}{\sqrt{2}}\in \left(-1, -\frac{1}{\sqrt{2}} \right) \cup \left( -\frac{1}{\sqrt{2}}, 1 \right)$. We may write
$$ \sqrt{\lambda^2 + 2\lambda -1} = \sqrt{-2(1-t^2)} = \sqrt{2} \sqrt{1-t^2} i, $$
where the last equality means that we fix the prescribed value of the square root (we could choose this value in our
previous considerations, as well). Then
$$ \sqrt{\lambda^2 + 2\lambda -1} = \sqrt{2} \sqrt{1-(\cos(\arccos t))^2} i  =
\sqrt{2} \sin(\arccos t) i; $$
$$ \lambda + 1 \pm \sqrt{\lambda^2 + 2\lambda -1} = \sqrt{2} \left( \cos(\arccos t)) \pm i \sin(\arccos t) \right) = \sqrt{2} e^{\pm i \arccos t}. $$
By~(\ref{f2_40}) and the last equalities we get
$$ p_n(\sqrt{2}t-1) =  \frac{1}{\sqrt{2}t + 1} \sin\left( \frac{3\pi}{4} n \right) + T_n(t) + $$
$$ +
\left(
\frac{t}{\sqrt{1-t^2} i} - \frac{1}{ (\sqrt{2}t + 1) \sqrt{2}  \sqrt{1-t^2} i}
\right)
i \sin(n\arccos t), $$
where $n\in\mathbb{Z}_+$, $t\in \left(-1, -\frac{1}{\sqrt{2}} \right) \cup \left( -\frac{1}{\sqrt{2}}, 1 \right)$.
Therefore, formula~(\ref{f2_20}) holds.
$\Box$

By formula~(\ref{f2_20}) or by the recurrent relation~(\ref{f2_10}) one can calculate that
$$ p_2(\lambda) = 2\lambda (\lambda + 2),\ 
p_3(\lambda) = \sqrt{2} \lambda (2\lambda^2 + 6\lambda + 3). $$
We see that two subsequent polynomials have a common root $0$. Thus, these polynomials are not orthogonal on
the real line. 

\noindent
It is an interesting open problem to describe the distribution of zeros for polynomials $p_n$ from~(\ref{f2_20}).
For example, we can conjecture that all zeros of $p_n$ are real.

\section{Orthogonality relations for the associated polynomials.}

Let an arbitrary Jacobi-type pencil $\Theta = \left(
J_3, J_5, \alpha, \beta
\right)$ be given. 
Set
\begin{equation}
\label{f3_40}
u_n := J_3 \vec e_n = a_{n-1} \vec e_{n-1} + b_n \vec e_n + a_n \vec e_{n+1},
\end{equation}
\begin{equation}
\label{f3_50}
w_n := J_5 \vec e_n = \gamma_{n-2} \vec e_{n-2} + \beta_{n-1} \vec e_{n-1} + \alpha_n \vec e_n + \beta_n \vec e_{n+1}
+ \gamma_n \vec e_{n+2},\qquad n\in\mathbb{Z}_+. 
\end{equation}
Here and in what follows by $\vec e_k$ with negative $k$ we mean (vector) zero.
Since $a_n > 0$, then vectors $\vec e_0$, $u_n$, $n\in\mathbb{Z}_+$ are linearly independent.
Moreover $\mathop{\rm Lin}\nolimits \{ \vec e_0, u_0, u_1, u_2, ... \} = l_{2,fin}$.
The following operator:   
$$ A f = \frac{\zeta}{\alpha} (\vec e_1 - \beta \vec e_0)
+ 
\sum_{n=0}^\infty \xi_n w_n, $$
\begin{equation}
\label{f3_60}
f = \zeta \vec e_0 + \sum_{n=0}^\infty \xi_n u_n\in l_{2,fin},\quad \zeta, \xi_n\in\mathbb{C}, 
\end{equation}
with $D(A) = l_{2,fin}$
is said to be \textbf{the associated operator for the Jacobi-type pencil $\Theta$}.
Notice that in the sums in~(\ref{f3_60}) only finite number of $\xi_n$ are nonzero. We shall always assume this in the case
of elements from the linear span.

Thus, the operator $A$ is linear and densely defined in $l_2$. In particular, we have:
\begin{equation}
\label{f3_70}
A u_n = w_n,\qquad n\in\mathbb{Z}_+, 
\end{equation}
\begin{equation}
\label{f3_80}
A \vec e_0 = \frac{1}{\alpha} \left( \vec e_1 - \beta \vec e_0 \right).
\end{equation}
By~(\ref{f3_70}),(\ref{f3_40}),(\ref{f3_50}) we obtain that
$$ \gamma_{n-2} \vec e_{n-2} + \beta_{n-1} \vec e_{n-1} - a_{n-1} A \vec e_{n-1} + \alpha_n \vec e_n - b_n A \vec e_n +
$$
\begin{equation}
\label{f3_90}
+ \beta_n \vec e_{n+1} - a_n A \vec e_{n+1} + \gamma_n \vec e_{n+2} = 0,\qquad n\in\mathbb{Z}_+.
\end{equation}
Consider the following equation:
$$ \gamma_{n-2} y_{n-2} + \beta_{n-1} y_{n-1} - a_{n-1} A y_{n-1} + \alpha_n y_n - b_n A y_n +
$$
\begin{equation}
\label{f3_100}
+ \beta_n y_{n+1} - a_n A y_{n+1} + \gamma_n y_{n+2} = 0,\qquad n\in\mathbb{Z}_+,
\end{equation}
with respect to unknown vectors $\{ y_n \}_{n=0}^\infty$, $y_n\in l_{2,fin}$, vectors $y_k$ with
negative $k$ are zero.
It is clear that the solution of~(\ref{f3_100}) is uniquely determined by $y_0, y_1$.
Vectors $\{ \vec e_n \}_{n=0}^\infty$ form a solution of~(\ref{f3_100}).

For an arbitrary non-zero polynomial $f(\lambda)\in\mathbb{P}$ of degree $d\in\mathbb{Z}_+$,
$f(\lambda) = \sum_{k=0}^d d_k \lambda^k$, $d_k\in\mathbb{C}$, we set
\begin{equation}
\label{f3_105}
f(A) = \sum_{k=0}^d d_k A^k.
\end{equation}
Here $A^0 := E|_{l_{2,fin}}$.
Since $A l_{2,fin} \subseteq l_{2,fin}$, then $D(f(A)) = l_{2,fin}$.
For $f(\lambda) \equiv 0$, we set
\begin{equation}
\label{f3_107}
f(A) = 0|_{l_{2,fin}}.
\end{equation}
The correspondence $f\mapsto f(A)$ is additive and multiplicative:
for arbitrary $f,g\in\mathbb{P}$ 
\begin{equation}
\label{f3_108}
(f+g)(A) = f(A) + g(A),\quad (fg)(A) = f(A)g(A).
\end{equation}
Rewrite~(\ref{f1_30}) with the operator argument $A$ and apply to $\vec e_0$. We obtain that
$\widetilde y_n = p_n(A) \vec e_0$ is a solution of~(\ref{f3_100}).
Notice that $\widetilde y_0 = \vec e_0$, 
$$ \widetilde y_1 = p_1(A) \vec e_0 = \vec e_1. $$
We conclude that solutions $\{ \vec e_n \}_{n=0}^\infty$, $\{ \widetilde y_n \}_{n=0}^\infty$ of~(\ref{f3_100})
coincide and
\begin{equation}
\label{f3_110}
\vec e_n = p_n(A) \vec e_0,\qquad n\in\mathbb{Z}_+.
\end{equation}
Therefore
\begin{equation}
\label{f3_120}
\left( p_n(A) \vec e_0, p_m(A) \vec e_0 \right)_{l_2} = \delta_{n,m},\qquad n,m\in\mathbb{Z}_+.
\end{equation}

\begin{center}
{\large\bf 
Orthogonal polynomials related to some Jacobi-type pencils.}
\end{center}
\begin{center}
{\bf S.M. Zagorodnyuk}
\end{center}

In this paper we study a generalization of the class of orthogonal polynomials on the real line.
These polynomials satisfy the following relation:
$(J_5 - \lambda J_3) \vec p(\lambda) = 0$,
where 
$J_3$ is a Jacobi matrix and
$J_5$ is a semi-infinite real symmetric five-diagonal matrix with positive numbers on the second subdiagonal,
$\vec p(\lambda) = (p_0(\lambda), p_1(\lambda), p_2(\lambda),\cdots)^T$, the superscript $T$ means the transposition,
with the initial conditions $p_0(\lambda) = 1$, $p_1(\lambda) = \alpha \lambda + \beta$, $\alpha > 0$, $\beta\in\mathbb{R}$.
Some orthonormality conditions for the polynomials $\{ p_n(\lambda) \}_{n=0}^\infty$ are obtained.
An explicit example of such polynomials is constructed.
}


\begin{thebibliography}{99}


\bibitem{cit_3000_Chihara}
Chihara, T. S. An introduction to orthogonal polynomials. Mathematics and its Applications, Vol. 13. 
Gordon and Breach Science Publishers, New York-London-Paris, 1978. xii+249 pp.

\bibitem{cit_3500_Ch_Z}
Choque Rivero, Abdon E.; Zagorodnyuk, Sergey M. Orthogonal polynomials on rays: Christoffel's formula. 
Bol. Soc. Mat. Mexicana (3) 15 (2009), no. 2, 149--164.

\bibitem{cit_4000_D_P_S}
Damanik, David; Pushnitski, Alexander; Simon, Barry. The analytic theory of matrix orthogonal polynomials. 
Surv. Approx. Theory 4 (2008), 1--85.

\bibitem{cit_5000_Ismail}
Ismail, Mourad E. H. Classical and quantum orthogonal polynomials in one variable. With two chapters by Walter Van Assche. 
With a foreword by Richard A. Askey. Encyclopedia of Mathematics and its Applications, 98. Cambridge University Press, Cambridge, 2005. xviii+706 pp.

\bibitem{cit_7000_Markus}
Markus, A. S. Introduction to the spectral theory of polynomial operator pencils. Translated from the Russian by H. H. McFaden. 
Translation edited by Ben Silver. With an appendix by M. V. Keldysh. Translations of Mathematical Monographs, 71. American Mathematical Society, 
Providence, RI, 1988. iv+250 pp.



\bibitem{cit_10000_Simon}
Simon, Barry. Szeg\"o's theorem and its descendants. Spectral theory for $L\sp 2$ perturbations of orthogonal polynomials. 
M. B. Porter Lectures. Princeton University Press, Princeton, NJ, 2011. xii+650 pp.

\bibitem{cit_20000_Suetin}
Suetin, P. K. Classical orthogonal polynomials.
Third edition. Fizmatlit, Moscow, 2005. 480 pp. (Russian)


\bibitem{cit_50000_Gabor_Szego}
Szeg\"o, G\'abor. Orthogonal polynomials. Fourth edition. 
American Mathematical Society, Colloquium Publications, Vol. XXIII. American Mathematical Society, Providence, R.I., 1975. xiii+432 pp.

\bibitem{cit_70000_T}
Teschl, Gerald. Jacobi operators and completely integrable nonlinear lattices. Mathematical Surveys and Monographs, 72. 
American Mathematical Society, Providence, RI, 2000. xvii+351 pp.

\bibitem{cit_90000_Z}
Zagorodnyuk, Sergey M. On generalized Jacobi matrices and orthogonal polynomials. New York J. Math. 9 (2003), 117--136 (electronic).


\end{thebibliography}
\end{document}